\documentclass[12pt]{article}
\usepackage{latexsym}
\usepackage{amsfonts}
\usepackage{amsmath}
\usepackage{amssymb}
\usepackage{amscd}
\usepackage{bbm}

\newcommand{\diam}{\mathop{\rm diam}}

\newcommand{\R}{\mathbb{R}}

\newcommand{\inr}[1]{\left<#1\right>}

\newcommand{\E}{\mathbb{E}}

\newcommand{\eps}{\varepsilon}
\newcommand{\ep}{\varepsilon}

\newcommand{\conv}{\mathop{\rm conv}}
\newcommand{\supp}{\mathop{\rm supp}}

\newcommand{\Pp}{\mathbb{P}}

\newcommand{\lstar}{\ell_*}

\makeatletter
\@addtoreset{equation}{section}

\newtheorem{Theorem}{Theorem}[section]
\newtheorem{Lemma}[Theorem]{Lemma}

\newtheorem{Definition}[Theorem]{Definition}

\newtheorem{Corollary}[Theorem]{Corollary}
\newtheorem{Remark}[Theorem]{Remark}

\numberwithin{equation}{section} 

\def \proof {\noindent {\bf Proof.}\ \ }

\def \endproof
{{\mbox{}\nolinebreak\hfill\rule{2mm}{2mm}\par\medbreak}}

\date{June 12, 2005}
\title{Reconstruction  and subgaussian operators}
 \author{
Shahar MENDELSON${}^{1}$ \and Alain PAJOR$^{2}$
       \and
Nicole TOMCZAK-JAEGERMANN${}^{3}$ }

\begin{document}
\maketitle
\footnotetext{{\it 2000 MSC-classification:}
    46B07,  41A45, 94B75, 52B05, 62G99}
\footnotetext[1]{Partially supported by an Australian Research
Council Discovery grant.}
 \footnotetext[2]{Partially supported by
an Australian Research Council grant.}
\footnotetext[3]{This author holds the Canada Research Chair in
Geometric Analysis.}

\noindent {\bf Abstract:}

We present a randomized method to approximate any vector $v$ from
some set $T \subset \R^n$. The data one is given is the set $T$ and
$k$ scalar products $(\inr{X_i,v})_{i=1}^k$, where $(X_i)_{i=1}^k$
are i.i.d. isotropic subgaussian random vectors in $\R^n$, and $k
\ll n$. We show that with high probability, any $y \in T$ for which
$(\inr{X_i,y})_{i=1}^k$ is close to the data vector
$(\inr{X_i,v})_{i=1}^k$ will be a good approximation of $v$, and
that the degree of approximation is determined by a natural
geometric parameter associated with the set $T$.

We also investigate a random method to identify exactly any vector
which has a relatively short support using linear subgaussian
measurements as above. It turns out that our analysis, when applied to
$\{-1,1\}$-valued vectors with i.i.d, symmetric entries, yields new
information on the geometry of faces of random $\{-1,1\}$-polytope; we
show that a $k$-dimensional random $\{-1,1\}$-polytope with $n$
vertices is $m$-neighborly for very large $m\le {c k/\log (c'
  n/k)}$.

The proofs are based on new estimates on the behavior of the
empirical process $\sup_{f \in F} \left|k^{-1}\sum_{i=1}^k f^2(X_i)
-\E f^2 \right|$ when $F$ is a subset of the $L_2$ sphere. The
estimates are given in terms of the $\gamma_2$ functional with
respect to the $\psi_2$ metric on $F$, and hold both in exponential
probability and in expectation.

\section{Introduction} \label{sec:intro}

The aim of this article is to investigate the linear ``approximate
reconstruction" problem in $\R^n$. In such a problem, one is given
a set $T \subset \R^n$ and the goal is to be able to approximate
any unknown $v \in T$ using random linear measurements. In other
words, one is given the set of values $(\inr{X_i,v})_{i=1}^k$,
where $X_1,...,X_k$ are independent random vectors in $\R^n$
selected according to some probability measure $\mu$. Using this
information (and the fact that the unknown vector $v$ belongs to
$T$) one has to produce, with very high probability with respect
to $\mu^k$,  some $t \in T$, such that the Euclidean norm $|t-v|
\leq \eps(k)$ for $\eps(k)$ as small as possible. Of course, the
random sampling method has to be ``universal" in some sense and
not tailored to a specific set $T$; and it is natural to expect that
the degree of approximation $\eps(k)$ depends on some geometric
parameter associated with $T$.

Questions of a similar flavor have been thoroughly studied in
nonparametric statistics and statistical learning theory (see, for
example, \cite{BBL} and \cite{M,MT} and references therein). This
particular problem has been addressed by several authors (see
\cite{CT1, CT2, CT3, RV} for the most recent contributions), in a
rather restricted context. First of all, the sets considered were
either the unit ball in $\ell_1^n$ or the unit balls in weak
$\ell_p^n$ spaces for $0<p<1$ - and the proofs of the approximation
estimates depended on the choice of those particular sets. Second, the
sampling process was done when $X_i$ were distributed according to the
Gaussian measure on $\R^n$ or in \cite{CT1} for Fourier ensemble.

In contrast, we present a method which is very
general. Our results hold for {\it any} set $T \subset \R^n$, and
the class of measures that could be used is broad; it contains all
probability measures on $\R^n$ which are isotropic and
subgaussian, that is, satisfy that for every $y \in \R^n$,
$\E|\inr{X,y}|^2=|y|^2$, and the random variables $\inr{X,y}$ are
subgaussian with constant $\alpha |y|$ for some $\alpha \geq 1$.
(see Definition \ref{def:measure}, below). This class of measures
contains, among others, the Gaussian measure on $\R^n$, the
uniform measure on the vertices of the combinatorial cube and the
normalized volume measure on various convex, symmetric bodies
(e.g. the unit balls of $\ell_p^n$ for $2 \leq p \leq \infty$).

It turns out that the key parameter in the estimate on the degree
of approximation $\eps(k)$ is indeed geometric in nature.
Moreover, the analysis of the approximation problem is centered
around the way the random operator $\Gamma=\sum_{i=1}^k
\inr{X_i,\cdot}e_i$ (where $(e_i)_{i=1}^k$ is the standard basis
in $\R^k$) acts on subsets of the unit sphere.

Our geometric approach, when applied to the sets $T$ considered in
\cite{CT1}, yields the optimal estimates on $\eps(k)$, and with
better probability estimates (of the order of $1-\exp(-ck)$), even
when the sampling is done according to an arbitrary isotropic,
subgaussian measure. Moreover, our result is more robust that the
one from \cite{CT1} in the following sense. The reconstruction
method suggested in \cite{CT1} is to find some $y \in T$ such that
for every $1 \leq i \leq k$, $\inr{X_i,v}=\inr{X_i,y}$, and then to
show that such  $v$
and $y$ must necessarily be   close. From our theorem it is clear
that one can choose any $y \in T$ for which $\sum_{i=1}^k
\inr{X_i,y-v}^2$ is relatively small, which is far more stable
algorithmically.

For the moment, let us present a simple version of the main
result we prove in this direction,  and to that end we require the
following notation. Let $(g_i)_{i=1}^n$ be independent standard
Gaussian random variables. Let $T\subset\R^n$ be a star-shaped
set (i.e. for every $t\in T$ and $0 \le \lambda\le 1$, $\lambda t\in T$)
and consider the following geometric parameter
$$
\lstar(T)=\E \sup_{t \in T}
   \left| \sum_{i=1}^n g_i t_i\right|
$$
which is, up to a factor of the order of $2\sqrt n$, the mean width of
the body $T$. We now define a more sensitive parameter, which as we
will see in this article, is the right parameter to control the error
term $\eps(k)$ for a star-shaped set $T$:
$$
r^*_k (\theta,T) := \inf \Bigl\{ \rho>0 : \rho \ge c \, \alpha^2 \lstar
(T\cap \rho S^{n-1})/  \theta\sqrt k  \Bigr\}.
$$

We may now state a version of our main result
concerning approximate reconstruction.

\medskip\noindent{\bf Theorem A }
{\it
There exist an absolute constant $c_1$ for which the following
holds. Let $T$ be a star-shaped subset of $\R^n$. Let $\mu$ be an
isotropic, subgaussian measure with constant $\alpha$ on $\R^n$ and
set $X_1,...,X_k$ be independent, selected according to $\mu$. Then,
with probability at least $1-\exp(-c_1k/\alpha^4)$, every $y,v \in
T$ satisfy that
$$
|y-v| \leq 2\left(\frac{1}{k} \sum_{i=1}^k
\left(\inr{X_i,v}-\inr{X_i,y}\right)^2 \right)^{1/2} +r_k^*(1/2, T-T).
$$
}

The parameter $r^*_k (\theta,T)$ can be estimated for unit ball of
classical normed or quasi-normed spaces (see Section
\ref{sec:app-rec}). In particular, if $T=B_1^n$ then with the same
hypothesis and probability as above, one has
$$
|y-v| \leq 2\left(\frac{1}{k} \sum_{i=1}^k
\left(\inr{X_i,v}-\inr{X_i,y}\right)^2 \right)^{1/2} +
c \alpha^2 \sqrt{\frac{\log({c\alpha^4 n/k})}{k}}
$$
where $c >0 $ is an absolute constant; this leads to the optimal
estimate for $\ep(k)$ for that set.

The main idea in the proof of this theorem is the fact that
excluding a set with exponentially small probability, the random
operator $\frac{1}{\sqrt{k}}\sum_{i=1}^k \inr{X_i,\cdot}$ is a very
good isomorphism on elements of $T$ whose Euclidean norm is large
enough (see Section \ref{random_proj} for more details).

A question of a similar nature we investigate here focuses on ``exact
reconstruction" of vectors in $\R^n$ that have a short support.
Suppose that $z \in \R^n$ is in the unit Euclidean ball, and has a
relatively short support $m <<n$. The aim is to use a random sampling
procedure to identify $z$ exactly, rather than just to approximate it.
The motivation for this problem comes from {\it error correcting
  codes}, in which one has to overcome random noise that corrupts a
part of a transmitted signal. The noise is modelled by adding to the
transmitted vector the noise vector $z$.  The assumption that the
noise does not change many bits in the original signal implies that
$z$ has a short support, and thus, in order to correct the code, one
has to identify the noise vector $z$ exactly. Since error correcting
codes are not the main focus of this article, we will not explore this
topic further, but rather refer the interested reader to \cite{MS,
  CT2, RV} and references therein for more information.

In the geometric context we are interested in, the problem has
been studied in \cite{ CT2, RV}, where it was shown that if $z$
has a short support relative to the dimension $n$ and the size of
the sample $k$, and if $\Gamma$ is a $k \times n$ matrix whose
entries are independent, standard Gaussian variables, then with
probability at least $1-\exp(-ck)$, the minimization problem
$$
(P)\qquad  \min \|v\|_{\ell_1^n} \ \ \ {\rm for} \ \ \Gamma v =
\Gamma z,
$$
has a unique solution, which is $z$. Thus, solution to this
minimization
problem will pin-point the ``noise vector" $z$.
The idea of using this minimization problem was first suggested in
\cite{CDS}.

We extend this result to a general random matrix whose rows are
$X_1,...,X_k$, sampled independently according to an isotropic,
subgaussian  measure.

\medskip\noindent{\bf Theorem B }
{\it
Let $\Gamma$ be as above.  With probability at least $1-
\exp(- c_1\,  k/\alpha ^4)$, any vector $z$ whose support has cardinality
 less than $ {c_2k/\log (c_3 n/k)}$ is the unique minimizer of the
problem $(P)$, where $c_1, c_2, c_3$ are absolute constants.
}

\medskip
Interestingly enough, the same analysis yields some information on
the geometry of $\{-1,1\}$ random polytopes. Indeed, consider the
$k \times n$ matrix $\Gamma$ whose entries are independent
symmetric $\{-1,1\}$ valued random variables. Thus, $\Gamma$ is a
random operator selected according to the uniform measure on the
combinatorial cube $\{-1,1\}^n$, which is an isotropic,
subgaussian measure with constant $\alpha=1$. The columns of
$\Gamma$, denoted by $v_1,...,v_n$ are vectors in $\{-1,1\}^k$ and
let $K^+ = {\rm conv}(v_1,...,v_n)$ be the convex
polytope generated by $\Gamma$; $K^+$ is thus called a random
$\{-1,1\}$-polytope.

A convex polytope is called $m$-neighborly if
any set of less than $m$ of its vertices
is the vertex set of a face
(see \cite{Z}). The following result yields the surprising fact
that a random $\{-1,1\}$-polytope  is $m$-neighborly for a
relatively large $m$. In particular, it will have the maximal
number of $r$-faces for $r\le m$.

\medskip\noindent{\bf Theorem C }
{\it
  There exist absolute constants $c_1,c_2, c_3$ for which the
  following holds. For $1\le k\le n$, with probability larger than $1-
  \exp(- c_1 k)$, a $k$-dimensional random $\{-1,1\}$ convex polytope
  with $n$ vertices is $m$-neighborly provided that
$$
m\le {c_2k\over \log (c_3 n/k)}\cdot
$$
}

\medskip

The main technical tool we use throughout this article,  which is
of independent interest, is an estimate of the behavior of the
supremum of the empirical process $f \to Z_f =k^{-1}\sum_{i=1}^k
f^2(X_i) - 1$ indexed by a subset of the $L_2$ sphere. That is,
$$
\sup_{f \in F} \left|\frac{1}{k} \sum_{i=1}^k f^2(X_i) -1 \right|,
$$
where $X_1,...,X_k$ are independent, distributed according to the
a probability measure $\mu$, and under the assumption that every
$f \in F$ satisfies that $\E f^2 =1$. We assume further that $F$
is a bounded set with respect to the $\psi_2$ norm, defined for a
random variable $Y$ by
$$
\|Y\|_{\psi_2} = \inf \left\{ u>0 : \E \exp(Y^2/u^2) \leq 2
\right\}.
$$

To formulate the following result, we require the notion of the
$\gamma_p$ functional \cite{Tal-book}.
For a metric space $(T,d)$, an {\it admissible sequence} of $T$ is a
collection of subsets of $T$, $\{T_s : s \geq 0\}$, such that for
every $s \geq 1$, $|T_s|=2^{2^s}$ and $|T_0|=1$.  for $p=1, 2$, we
define the $\gamma_p$ functional by
$$
\gamma_p(T,d) =\inf \sup_{t \in T} \sum_{s=0}^\infty
2^{s/p}d(t,T_s),
$$
where the infimum is taken with respect to all admissible sequences of
$T$.

\medskip\noindent{\bf Theorem D }
{\it
There exist absolute constants $c_1, c_2, c_3$ and for which the
following holds. Let $(\Omega,\mu)$ be a probability space, set
$F$ be a subset of the unit sphere of $L_2(\mu)$ and assume that
${\rm diam}(F, \| \ \|_{\psi_2}) = \alpha$. Then, for any
$\theta>0$ and $k\ge 1$ satisfying
$$
c_1\alpha \gamma_2(F, \| \
\|_{\psi_2})\leq \theta\sqrt{k},
$$
with probability at least $1-\exp(-c_2\theta^2k/\alpha^4)$,
$$
\sup_{f \in F} |Z_f| \leq \theta.
$$
Moreover, if $F$ is symmetric, then
$$
\E\sup_{f \in F} |Z_f| \leq c_3\alpha^2 \frac{\gamma_2(F, \| \
\|_{\psi_2})}{\sqrt{k}}.
$$
}

Theorem D 
improves a result of a similar
flavor from \cite{KM} in two ways. First of all, the bound on the
probability is exponential in the sample size $k$ which was not
the case in \cite{KM}. Second, we were able to bound the expectation
of the supremum of the empirical process using only a $\gamma_2$
term. This fact is surprising because the expectation of the
supremum of empirical processes is usually controlled by two
terms; the first one bounds the subgaussian part of the process
and is captured by the $\gamma_2$ functional with respect to the
underlying metric. The other term is needed to control
sub-exponential part of the empirical process, and is bounded by
the $\gamma_1$ functional with respect to an appropriate metric
(see \cite{Tal-book} for more information on the connections
between the $\gamma_p$ functionals and empirical processes).
Theorem  D
shows that the expectation of the
supremum of $|Z_f|$ behaves as if $\{Z_f\,:\, f\in F\}$ were a
subgaussian process with respect to the $\psi_2$ metric (although
it is not), and this is due to the key fact that all the functions
in $F$ have the same second moment.

We end this introduction with the organization of the article. In
Section \ref{empirical} we present the proof of Theorem D
and some of its corollaries we require. In
Section \ref{random_proj} we illustrate these results in the case
of linear processes which corresponds to linear measurements. In
Section \ref{sec:app-rec} we investigate the approximate
reconstruction problem for a general set, and in Section
\ref{sec:exact} we present a proof for the exact reconstruction
scheme, with its application to the geometry of random
$\{-1,1\}$-polytopes.

Throughout this article we will use letters such as $c,c_1,..$ to
denote absolute constants which may change depending on the context.
We denote by $(e_i)_{i=1}^n$ the canonical basis of $\R^n$, by $|x|$
the Euclidean norm of a vector $x$ and by $B_2^n$ the Euclidean unit
ball. We also denote by $|I|$ the cardinality of a finite set $I$.

\medskip
\noindent{\it Acknowledgement:}
A part of this work was done when the second and the third authors
were visiting the Australian National University, Canberra; and when
the first and the third authors were visiting Universit\'{e} de
Marne-la-Vall\'{e}e.  They wish to thank these institutions for their
hospitality.

\section{Empirical Processes}
\label{empirical}

In this section we present some results in empirical processes that
will be central to our analysis. All the results focus on
understanding the process
$ Z_f=\frac{1}{k}\sum_{i=1}^k f^2(X_i) -\E f^2 $,
where $k \ge 1$ and $X_1,...,X_k$ are independent random variables
distributed according to a probability measure $\mu$. In particular,
we investigate the behavior of $\sup_{f \in F} |Z_f|$ in terms of
various metric structures on $F$, and under the key assumption that
every $f \in F$ has the same second moment. The parameters involved
are standard in {\it generic chaining} type estimates (see
\cite{Tal-book} for a comprehensive study of this topic).

\smallskip

Recall that the $\psi_p$ norm of a random variable $X$ is
defined as
$$
\|X\|_{\psi_p}=\inf \left\{u >0:
  \E\exp\left({|X|^p}/{u^p}\right) \leq 2 \right\}.
$$
It is standard to verify (see for example \cite{VW}) that if
$X$ has a bounded $\psi_2$ norm, then it is subgaussian with parameter
$c\|X\|_{\psi_2}$ for some absolute constant $c$.
More generally, a bounded $\psi_p$ norm implies that
$X$ has a tail behavior, $\Pp (|X| >u)$,
of the type $\exp (- cu^p/\|X\|_{\psi_p})$.

Our first fundamental ingredient is the well known Bernstein's
inequality which we shall use in the form of  a $\psi_1$ estimates
(\cite{VW}).

\begin{Lemma}
    \label{Bernst-psi1}
Let $Y_1,...,Y_k$ be independent random variables with zero mean
such that for some $b >0$ and every $i$, $\|Y_i\|_{\psi_1} \le b$.
Then, for any  $u>0$,
\begin{equation} \label{eq:Bern}
\Pp \Bigl(\bigl|\frac{1}{k}\sum_{i=i}^k Y_i \bigr| >u
\Bigr) \leq 2\exp \left(-c\, k\min \Bigl(
\frac{u}{b},\, \frac{u^2}{b^2}\Bigr)\right),
\end{equation}
where $c >0$  is an absolute constant.
\end{Lemma}

We will be interested in classes of functions $F \subset L_2(\mu)$
bounded in the $\psi_2$-norm; we assume without loss of generality
that $F$ is symmetric and we let $\diam(F, \|\ \|_{\psi_2}) : = 2
\sup_{f \in F} \|f\|_{\psi_2}$.  Additionally, in many technical
arguments we shall often consider classes $F \subset S_{L_2}$, where
$S_{L_2} = \{ f : \|f\|_{L_2} =1 \}$ is the unit sphere in $L_2(\mu)$.

\medskip

Let $X_1, X_2, \ldots$ be independent
random variables distributed according to $\mu$.  Fix  $k \ge 1$ and
for $f \in F$ define the random variables  $Z_f$ and  $W_f$ by
$$
 Z_f=\frac{1}{k}\sum_{i=1}^k f^2(X_i)- \E f^2
\qquad \mbox{\rm and} \qquad
W_f=\left(\frac{1}{k} \sum_{i=1}^k f^2(X_i)\right)^{1/2}.
$$

\smallskip

The first lemma follows easily from Bernstein's inequality.  We state it
in the form convenient for further use, and give a proof of one
part, for completeness.

\begin{Lemma}
  \label{lemma:psi-2-est}
  There exists an absolute constant $c_1 >0$ for which the following
  holds.  Let $F \subset S_{L_2}$,
  $\alpha = \diam (F, \|\ \|_{\psi_2})$ and set $k \ge 1$.  For every
  $f, g \in F$ and every $u \geq 2$ we have
$$
\Pp\bigl(W_{f-g} \geq  u\|f-g\|_{\psi_2} \bigr) \le
 2 \exp(-c_1 k u ^2).
$$
Also, for every $u >0$,
$$
     \Pp \left( \left|Z_f - Z_g\right|  \ge  u\, \alpha\,
        \|f-g\|_{\psi_2} \right)
   \le 2 \exp \bigl(-c_1 k \min (u, u^2)\bigr),
$$
and
$$
     \Pp \left( \left|Z_f \right|  \ge  u\, \alpha^2 \right)
   \le 2 \exp \bigl(-c_1 k \min (u, u^2)\bigr).
$$
\end{Lemma}

\proof We show the standard proof of the first estimate.  Other
estimates are proved similarly (see e.g., \cite{KM}, Lemma 3.2).

Clearly,
$$
\E \,W_{f-g}^2 = \frac{1}{k} \E \sum_{i=1}^k   (f-g)^2 (X_i)
= \E (f-g)^2 (X_1) = \|f-g\|_{L_2}^2.
$$
Applying Bernstein's inequality  it follows that for $t >0$,
\begin{eqnarray*}
  \lefteqn{\Pp \left(\left|W_{f-g}^2 - \|f-g\|_{L_2}^2 \right| \ge t
    \right)}\\
&\le& 2\, \exp \left(- c k \min \left(\frac{t}{\| (f-g)^2\|_{\psi_1}},
\Bigl(\frac{t}{\| (f-g)^2 \|_{\psi_1}}\Bigr)^2\right)\right).
\end{eqnarray*}
Since $\|h^2 \|_{\psi_1} = \|h \|_{\psi_2}^2 $ for every function $h$,
then
letting $t = (u^2-1) \| f-g\|_{\psi_2}^2$,
\begin{eqnarray*}
  \Pp \Bigl( W_{f-g}^2 &\ge & u^2 \| f-g\|_{\psi_2}^2 \Bigr)\\
& \le  &
\Pp \Bigl( W_{f-g}^2 - \|f-g\|_{L_2}^2  \ge
 (u^2-1) \| f-g\|_{\psi_2}^2
 \Bigr) \\
& \le & 2\, \exp \left(- c k \min (u^2/2, u^4/4) \right),
\end{eqnarray*}
as promised.
\endproof

\medskip

Now we return to one of the basic notions used in this paper, that
of the $\gamma_2$-functional.
Let $(T,d)$ be a metric space. Recall that an {\it admissible sequence} of
$T$ is a collection of subsets of $T$, $\{T_s : s \geq 0\}$, such
that for every $s \geq 1$, $|T_s|=2^{2^s}$ and $|T_0|=1$.

\begin{Definition}
    \label{def:gamma-2}
For a metric space $(T,d)$  and $p = 1, 2$, define
$$
\gamma_p(T,d) =\inf \sup_{t \in T} \sum_{s=0}^\infty
2^{s/p}d(t,T_s),
$$
where the infimum is taken with respect to all admissible sequences of $T$.
In cases where the metric is clear from the context, we will denote
the $\gamma_p$ functional by $\gamma_p(T)$.
\end{Definition}

Set $\pi_s:T \to T_s$ to be a metric projection function onto $T_s$,
that is, for every $t \in T$, $\pi_s(t)$ is a nearest element to $t$
in $T_s$ with respect to the metric $d$. It is easy to verify by the
triangle inequality that for every admissible sequence and every $t \in T$,
$\sum_{s=0}^\infty 2^{s/2}d(\pi_{s+1}(t),\pi_{s}(t)) \leq ( 1 +
1/\sqrt 2) \sum_{s=0}^\infty 2^{s/2}d(t,T_s)$ and that $\diam(T,d) \leq
2\gamma_2(T,d)$.

\medskip

We say that a set $F$ is {\it star-shaped} if the fact that $f \in F$
implies that $\lambda f \in F$ for every $0 \leq \lambda \leq 1$.

The next Theorem shows that excluding a set with exponentially
small probability, $W_f$ is close to being an isometry
in the $L_2(\mu)$ sense for functions in $F$ that have a
relatively large norm.

\begin{Theorem}
  \label{thm:general}
  There exist absolute constants $c, \bar{c} >0$ for which the
  following holds.  Let $F \subset L_2(\mu)$ be star-shaped, $\alpha =
  \diam (F, \|\ \|_{\psi_2})$ and $k \ge 1$. For any $0 < \theta <1$,
  with probability at least $1- \exp(- \bar{c}\, \theta^2 k/\alpha
  ^4)$, then for all $f \in F$ satisfying $\E f^2 \ge
  r_k^*(\theta)^2$, we have
\begin{equation}
  \label{two_sided_gen}
   (1 - \theta) \E f^2 \le \frac{ W_f^2 }{ k} \le (1+\theta) \E f^2,
\end{equation}
where
\begin{equation}
  \label{rho_star0}
  r^*_k(\theta) = r^*_k (\theta, F) := \inf \Bigl\{ \rho>0 : \rho \ge
  c \, \alpha
\frac{\gamma_2 (F \cap \rho S_{L_2}, \|\ \|_{\psi_2})}
        { \theta \sqrt k } \Bigr\}.
\end{equation}
\end{Theorem}

The two-sided inequality (\ref{two_sided_gen}) is intimately related
to an estimate on $\sup_{f \in F} |Z_f|$, which, in turn, is based on
two ingredients. The first one shows, in the language of the standard
chaining approach, that one can control the ``end parts" of all
chains. its proof  is essentially the same as Lemma 2.3 from \cite{KM}.

\begin{Lemma}
  \label{lemma:chain}
  There exists an absolute constant $C$ for which the following holds.
  Let $F \subset S_{L_2}$,
  $\alpha = \diam (F, \|\ \|_{\psi_2})$ and $k \ge 1$.  There is
  $F' \subset F$ such that $|F'| \leq 4^k$ and with probability at
  least $1-\exp(- k)$, we have, for every $f \in F$,
  \begin{equation}
    \label{chain}
W_{f-\pi_{F'}(f)} \leq C {\gamma_2(F,  \| \ \|_{\psi_2})}/{\sqrt{k}},
  \end{equation}
where $\pi_{F'}(f)$ is a nearest point to $f$ in $F'$
with respect to the $\psi_2$ metric.
\end{Lemma}

\proof Let $\{F_s: s \geq 0\}$ be an ``almost optimal'' admissible
sequence of $F$.  Then for every $f \in F$,
$$
\sum_{s=0}^\infty
2^{s/2}\, \|\pi_{s+1}(f)-\pi_{s}(f)\|_{\psi_2} \leq 2 \gamma_2(F,  \| \
\|_{\psi_2} ).
$$

Let $s_0$ be the minimal integer such that $2^{s_0} > k$,
and let $F'=F_{s_0}$. Then  $|F'| \leq 2^{2k} = 4^k$.
Write
$$
f-\pi_{s_0}(f)=\sum_{s=s_0}^\infty
\left(\pi_{s+1}(f)-\pi_s(f)\right).
$$
Since $W$ is  sub-additive  then
$$
W_{f-\pi_{s_0}}(f) \leq \sum_{s=s_0}^\infty
W_{\pi_{s+1}(f)-\pi_s(f)}.
$$
For any $f \in F$, $s \ge s_0$ and $\xi \ge 2$, noting that $2^s>k$,
it follows  by  Lemma~\ref{lemma:psi-2-est}  that
\begin{equation}
  \label{term_est}
\Pp \left(W_{\pi_{s+1}(f)-\pi_s(f)} \geq \xi
\frac{2^{s/2}}{\sqrt{k}}\|\pi_{s+1}(f)-\pi_s(f)\|_{\psi_2}\right)
 \leq 2 \exp(-c_1 \xi^2 2^{s}).
\end{equation}

Since $|F_s| \leq 2^{2^s}$, there are at most
$2^{2^{s+2}}$  pairs of $\pi_{s+1}(f)$ and $\pi_s(f)$. Thus,
for every $s \ge s_0$,
the probability of the event from
(\ref{term_est})  holding for
some $f \in F$ is less than or equal  to $2^{2^{s+2}}\cdot
2 \exp (- c_1 \xi^2 2^{s}) \le \exp (2^{s+3} - c \xi^2 2^{s})$,
which, for $\xi \ge \xi_0 := \max(4/\sqrt {c_1}, 2)$, does not exceed
$\exp (- c_1 \xi^2 2^{s-1})$.

Combining these estimates together it follows that
$$
W_{f-\pi_{s_0}(f)} \leq \xi \sum_{s=s_0}^\infty
\frac{2^{s/2}}{\sqrt{k}}\|\pi_{s+1}(f)-\pi_s(f)\|_{\psi_2} \le
2 \xi\, \frac{\gamma_2\bigl(F, \|\ \|_{\psi_2}\bigr)}{\sqrt{k}},
$$
outside a set of probability
$$
\sum_{s=s_0}^\infty \exp (- c_1 \xi^2 2^{s-1}) \le \exp(- c_1
\xi^2  2^{s_0} /4).
$$
We complete the proof setting, for example $\xi =
\max(\xi_0, 2/\sqrt{c_1})$ and recalling that $2^{s_0} > k$.
\endproof

\begin{Remark}
  \label{rem:expect}
{\rm The proof of the lemma shows that there exist absolute constants
  $c', c'' >0$ such that for every $\xi \ge c'$,
$$
\Pp \Biggl(\sup_{f \in F} W_{f-\pi(f)} \geq \xi
\frac{\gamma_2\bigl(F, \|\ \|_{\psi_2}\bigr)}{\sqrt{k}} \Biggr)
 \leq \exp(-c'' \xi^2 k),
$$
a fact which will be used later.}
\end{Remark}

\medskip

The next lemma estimates the supremum $\sup_{f \in F'} |Z_f|$, where
the supremum is taken over a subset $F'$ of $F$ of a relatively small
cardinality, or in other words, over the ``beginning part'' of a
chain.  However, in order to get an exponential in $k$ estimates
on probability, we require a separate argument (generic chaining) for
the ``middle part'' of a chain while for the ``very beginning'' it is
sufficient to use a standard concentration estimate.

\begin{Lemma}
   \label{lemma:small-set}
   There exist absolute constants $C$ and $c''' >0$ for which the
   following holds. Let $F \subset S_{L_2}$   and  $\alpha = \diam (F, \|\
   \|_{\psi_2})$.
   Let $k \ge 1$ and $F' \subset F$ such that $|F'| \leq 4^k$.  Then
   for every $ w > 0 $,
$$
\sup_{f \in F'}|Z_f| \leq C \alpha
  \frac{\gamma_2(F,  \| \ \|_{\psi_2})}{\sqrt{k}}+
\alpha^2 w,
$$
with probability larger than or equal to
$1 - 3 \exp(- c''' \,k\, \min (w, w^2) )$.
\end{Lemma}

\proof
Let $(F_s)_{s=0}^\infty$ be an almost optimal admissible sequence of
$F'$, set $s_0$ to be the minimal integer such that $2^{s_0} > 2 k$
and fix $s_1 \le  s_0$ to be determined later.
Since $|F'| \leq 4^k$, it follows that $F_s=F'$ for every $s \geq
s_0$,  and that
$$
Z_f-Z_{\pi_{s_1}(f)}=\sum_{s=s_1+1}^{s_0}
\left(Z_{\pi_s(f)}-Z_{\pi_{s-1}(f)}\right).
$$

By Lemma~\ref{lemma:psi-2-est}, for every $f \in F'$, $1 \le s
\leq s_0$ and $u >0$,
\begin{align*}
   \Pp \Bigl(  \left |Z_{\pi_{s}(f)}-Z_{\pi_{s-1}(f)}\right| &\geq
 u\alpha \sqrt{\frac{2^{s}}{k}}\,
\|\pi_{s+1}(x)-\pi_s(x)\|_{\psi_2}\Bigr) \\
& \leq   2\exp(-c_1 \min\bigl( (u \sqrt{2^s/k}), (u \sqrt{2^s/k})^2\bigr))\\
& \leq  2\exp(-c_1 \min(u,u^2)2^{s-2}).
\end{align*}
  (For the latter inequality observe that if  $s \leq s_0$
then $2^s/k \leq 4$, and thus  $\min\bigl( (u \sqrt{2^s/k}), (u
\sqrt{2^s/k})^2 \bigr) \ge \min (u, u^2 )\, 2^s/(4\,k)$.)

\smallskip

Taking $u$ large enough (for
example, $u = 2^5/c_1$ will suffice) we may ensure that
$$
 \sum_{s=s_1+1}^{s_0} 2^{2^{s+2}}\exp(-c_1 u2^{s-2})
\leq \sum_{s=s_1+1}^{s_0} \exp(-2^{s+3})
\leq
\exp(- 2^{s_1}).
$$
Therefore,
since there are at most
$2^{2^{s+2}}$  possible  pairs of $\pi_{s+1}(f)$ and $\pi_s(f)$,
there is a set of probability at most
$\exp(- 2^{s_1})$  such that outside  this set we
have
$$
\sup_{f \in F'}|Z_f-Z_{\pi_{s_1}(f)}|
\leq \frac{\alpha \, u }{\sqrt k}\,
\sum_{s=s_1}^{s_0} 2^{s/2}
\|\pi_{s+1}(x)-\pi_s(x)\|_{\psi_2} \le
c' \alpha \frac{\gamma_2(F)}{\sqrt{k}}.
$$

Denote $F_{s_1}$ by $F''$ and observe that $|F''| \le 2^{2^{s_1}}$.
Thus  the later estimate implies
$$
\sup_{f \in F'}|Z_f| \leq
c' \alpha \frac{\gamma_2(F)}{\sqrt{k}} +
\sup_{g \in F''}|Z_g|.
$$
Applying Lemma~\ref{lemma:psi-2-est}, for every $ w >0$ we get
$$
\Pp\left(|Z_{g}| \geq \alpha ^2 w \right) \leq
2\exp(- c_1   k \min(w, w^2)).
$$
Thus, given  $w >0$,
choose $s_1\le s_0$ to be the largest
integer such that  $2^{s_1} <  c_1 k \min(w, w^2)/2$.
Therefore,  outside  a set of probability  less than or equal to
$|F''|\, 2\exp(- c_1 k \min(w, w^2)) \le
\exp(- c_1 k \min(w, w^2)/2)$ we have
$|Z_g| \le \alpha ^2 w $ for all $g \in F''$.
To conclude, outside  a set of probability
$3\exp(- c_1 k \min(w, w^2)/2 )$,
$$
\sup_{f \in F'}|Z_f| \leq
c' \alpha \frac{\gamma_2(F)}{\sqrt{k}} + \alpha ^2 w,
$$
as required.
\endproof

\noindent {\bf Proof of Theorem~\ref{thm:general}.}
 Fix an arbitrary $\rho >0$, and for the purpose of this proof we let
 $ F(\rho) = F/\rho \cap S_{L_2}$, where $F/\rho=\{f/\rho : f\in F\}$.

Our first and main aim is to estimate $\sup_{f \in F(\rho)} |Z_f|$ on a
set of probability close to 1.

\smallskip

Fix $ u, w >0$ to be determined later.
Let $F' \subset F(\rho) $ be as Lemma~\ref{lemma:chain}, with $|F'|
\le 4^ k$. For every $f \in F(\rho) $ denote by $\pi(f)=\pi_{F'}(f)$ a
closest point to $f$ with respect to the $\psi_2$ metric on $F(\rho)$.
By writing $ f =( f -\pi(f))+\pi(f)$, it is evident that
$$
|Z_f| \leq W_{f-\pi(f)}^2+2W_{f-\pi(f)}W_{\pi(f)} + |Z_{\pi(f)}|,
$$
and  thus,
\begin{equation}
  \label{three_terms}
\sup_{f \in F(\rho)} |Z_f| \leq \sup_{f \in F(\rho)} W_{f-\pi(f)}^2
    + 2\sup_{f \in F(\rho)} W_{f-\pi(f)}  \sup_{g \in F'} W_{g}
     + \sup_{g \in F'}|Z_g|.
\end{equation}

Applying  Lemma~\ref{lemma:chain},
the first term in (\ref{three_terms}) is estimated using the fact that
$$
\sup_{f \in F(\rho)} W_{f-\pi(f)}
\le C \frac{\gamma_2(F(\rho), \| \ \|_{\psi_2})}{\sqrt k},
$$
with probability at least $1 - \exp(-  k)$.

\medskip

For every $f \in F(\rho) $ and every $u >0$ we have
$$
\Pp \left\{  W_f  \ge 1 + u \alpha ^2 \right\}  \le
\Pp \left\{  W_f^2  \ge    1 + u \alpha ^2 \right\}
\le \Pp \left\{ | Z_f|  \ge  u \alpha ^2 \right\}
$$
and, by Lemma~\ref{lemma:psi-2-est}, the latter probability  is at most
$ 2 \exp (-c\, k \min(u, u^2))$, where $c >0$ is an
absolute constant.

Combining  these two  estimates  with Lemma~\ref{lemma:small-set}
 and substituting into   (\ref{three_terms}),
$\sup_{f \in F(\rho)} |Z_f|$  is upper bounded by

\begin{eqnarray}
  \label{three_terms_est}
\lefteqn{
C^2 \frac{\gamma_2(F(\rho), \| \ \|_{\psi_2})^2} { k}
+  C \frac{\gamma_2(F(\rho), \| \ \|_{\psi_2})}{\sqrt k}
          (1+ u \alpha^2)} \nonumber\\
& & +   C'' \alpha
 \frac{\gamma_2(F(\rho), \| \ \|_{\psi_2})}{\sqrt k}   +
 \alpha^2 w,
\end{eqnarray}
with probability at least $1 - 2 e^{- k} - 2 e^{-c k
\min(u, u^2)} - 3 e^{-c k \min(w, w^2)}  $.

\smallskip

Given $0 < \theta < 1$ we want the condition $\sup_{f \in F(\rho)}
|Z_f| \le \theta $ to be satisfied with probability close to 1.
This can  be achieved  by imposing suitable conditions on the
parameters involved.  Namely, if  $u = 1/\alpha^2 < 1$
and if   $\rho > 0 $ and $  w >0$  satisfy
\begin{equation}
  \label{condition1}
  \tilde{C}\, \alpha\, \frac{\gamma_2(F(\rho), \| \ \|_{\psi_2})}{\sqrt k}
\le   \theta/4,
\qquad
C''' \alpha^2 w \le  \theta/4,
\end{equation}
where $\tilde {C} = \max(2\, C, C'')$, then
each of the last three terms in (\ref{three_terms_est})
is less than or equal
to $\theta /4$, and the first term is less than or equal to
$(\theta/4)^2$.

In order to ensure that (\ref{condition1}) holds, we let $w =
\min\bigl(1, \theta/(4 C''' \alpha^2)\bigr)$. The above discussion
shows that as long as $\rho$  satisfies
\begin{equation}
  \label{condition-rho}
  4 \,\tilde{C}\, \alpha\, \frac{\gamma_2(F(\rho), \| \ \|_{\psi_2})}
{\theta \sqrt k}  \le  1,
\end{equation}
then $\sup_{f \in F(\rho)} |Z_f| \le \theta$ on a set of measure larger
than or equal to $1 - 7 e^{- c k \theta^2/\alpha^4}$, where $c >0$ is
an absolute constant.
Hence, whenever $\rho$ satisfies (\ref{condition-rho}) then
(\ref{two_sided_gen}) holds for all $f \in F(\rho)$. Finally, note that
$\gamma_2(F(\rho), \|\ \|_{\psi_2}) = (1/\rho)\gamma_2(F \cap \rho
S_{L_2}, \|\ \|_{\psi_2})$, and  thus (\ref{condition-rho}) is equivalent
to the inequality in the definition of $ r^*_k(\theta) $.

\medskip

To conclude the proof, for a fixed  $0 < \theta <1$
set $r = r_k^*(\theta)$, with  $c = 4 \tilde{C}$ being
the constant from (\ref{condition-rho}).
Note that if $X_1,...,X_k$ satisfy (\ref{two_sided_gen}) for all $f \in
F(r)$ then, since $F$ is star-shaped, the homogeneity of this condition
implies that the same holds for all $f \in F$ with $\E f^2 \ge r^2$,
as claimed.
\endproof

\smallskip

Let us note two consequences
for the supremum of the process $Z_f$, which is of independent
interest.

\begin{Corollary}
  \label{est_sup_Z}
  There exist absolute constants $C', c' >0$ for which the
  following holds.  Let $F \subset S_{L_2}$, $\alpha = \diam (F,
  \|\ \|_{\psi_2})$ and $k \ge 1$.
With probability at least $1- \exp\left(- c'
    \gamma_2^2(F, \| \ \|_{\psi_2}) /\alpha^3\right)$ one has
$$
\sup_{f \in F} |Z_f| \le
C'\, \alpha\, \max\Bigl(\frac{\gamma_2(F, \| \ \|_{\psi_2})}{\sqrt k},
\frac{\gamma_2^2(F, \| \ \|_{\psi_2})}{  k}\Bigr).
$$
Moreover, if $F$ is symmetric,
$$
\E \sup_{f \in F} |Z_f| \le
C'\, \alpha\, \max\Bigl(\frac{\gamma_2(F, \| \ \|_{\psi_2})}{\sqrt k},
\frac{\gamma_2^2(F, \| \ \|_{\psi_2})}{  k}\Bigr).
$$
\end{Corollary}

\proof This follows from the proof of Theorem~\ref{thm:general} with
$\rho =1$.  More precisely, the first part is a direct consequence of
(\ref{three_terms_est}).

For the ``moreover part''
first use (\ref{three_terms}) for expectations, estimate the middle
term by Cauchy-Schwarz inequality and note that $W_g^2 \le 1 + Z_g$
for all $g \in F'$ to yield that in order to estimate $\E \sup_{f \in
  F} |Z_f|$ it suffice to bound
$$
\E \sup_{f \in F'} |Z_f|
\quad {\rm and} \quad
\E \sup_{f \in F} W_{f -\pi(f)}^2.
$$

For simplicity denote $\gamma_2(F, \| \ \|_{\psi_2})$ by $\gamma_2(F)$
and let us begin with the second term. Applying Remark~\ref{rem:expect}
and setting $G=\{f - \pi(f) : f \in F\}$ and $u = c'
\gamma_2(F)/\sqrt{k}$, where $c' $ is the constant from the remark, we
obtain
\begin{align*}
\int_0^\infty \Pp \left(\sup_{g \in G} W_g^2 \geq t\right)dt & \leq
u^2 + \int_{u^2}^\infty \Pp \left(\sup_{g \in G} W_g^2 \geq t
\right)dt
\\
& \leq u^2+u^2 \int_1^\infty \exp(-c'' v k). dv,
\end{align*}
where   the last inequality follows by changing the  integration
variable to $t=u^2 v$. This implies that
$$
\E \sup_{f \in F} W_{f-\pi(f)}^2 \leq C' \frac{\gamma_2^2(F,
\|\ \|_{\psi_2})}{k},
$$
for some absolute constant $C' >1$.

Next, we have to bound $\E \sup_{f \in F'} |Z_f|$, and to that end
we use
Lemma~\ref{lemma:small-set}.
Setting $u = 2 C \alpha \gamma_2(F)/\sqrt{k}$, and then changing
the integration variable to $t = u/2+\alpha^2w$, it is evident that
\begin{align*}
\int_0^\infty \Pp \left(\sup_{f \in F'} |Z_f| \geq t \right)dt
\leq & \, u + \int_{u}^\infty \Pp \left(\sup_{f \in F'}
       |Z_f| \geq t \right)dt \\
= & \,u + \alpha^2 \int_{u/2\alpha^2}^\infty \Pp \left(\sup_{f \in F'}
|Z_f|
\geq  \frac{u}{2} +\alpha^2 w \right)dw \\
\leq &\,  u + 3 \, \alpha^2
\int_{u/2\alpha^2}^\infty \exp\left(-c'''k\min(w^2,w)\right)dw.
\end{align*}
Changing variables in the last integral $w=ru/2\alpha^2$ and using the
fact that $\gamma_2(F) \geq 1$ for a symmetric set $F$ in the unit
sphere, the last
expression  is bounded above by
$$
u + (3/2) u
\int_{1}^\infty \exp\left(-c'''\min\left(\frac{r}{2\alpha^2},
\Bigl(\frac{r}{2\alpha^2}\Bigr)^2\right)\right)dr =
 C' u,
$$
where $C >0$ is an absolute constant.
\endproof

\section{Subgaussian Operators}
\label{random_proj}

We now illustrate the general result of Section~\ref{empirical} in the
case of linear processes, which was the topic that motivated our
study.  The processes corresponds then to random matrices with rows
distributed according to measures on $\R^n$ satisfying some natural
geometric conditions. Our result imply concentration estimates for
related random subgaussian operators, which eventually lead to the
desired reconstruction results for linear measurements for general
sets.

The fundamental result that allows us to pass from the purely metric
statement of the previous section to the geometric result we present
below follows from Talagrand's lower bound on the expectation of the
supremum of a Gaussian process in terms of $\gamma_2$ of the indexing
set. To present our result, the starting point is the fundamental
definition of the $\lstar$-functional (which is in fact the so-called
$\ell$-functional of a polar set).

\begin{Definition}
   \label{def:gauss}
   Let $T \subset \R^n$ and let $g_1,...,g_n$ be independent  standard
   Gaussian random variables.  Denote by $\lstar(T)=\E \sup_{t \in T}
   \left| \sum_{i=1}^n g_i t_i\right|$, where $ t = (t_i)_{i=1}^n \in
   \R^n$.
\end{Definition}

There is a close connection between the $\lstar$- and $\gamma_2$-
functionals given by the majorizing measure Theorem.
Let $\{G_t: t \in T\}$ be a Gaussian process indexed by a set $T$, and
for every $s, t \in T$, let $d^2 (s,t) = \E|G_s-G_t|^2$.
Then
$$
c_2 \gamma_2(T,d) \leq \E \sup_{t \in T} |G_t| \leq c_3 \gamma_2(T,d),
$$
where $c_2, c_3>0$ are  absolute constants.
The upper bound is due to Fernique \cite{F} and the lower bound was
established by Talagrand \cite{Ta}. The proof of both parts and the
most recent survey on the topic can be found in \cite{Tal-book}

In particular, if  $T \subset \R^n$ and  $G_t = \sum g_i t_i$, then $d(s, t)
= |s-t|$, and thus
\begin{equation}
  \label{maj_meas}
  c_2 \gamma_2(T,|\ |) \leq  \lstar (T)\leq c_3 \gamma_2(T,|\ |),
\end{equation}

\begin{Definition}
   \label{def:measure}
   A probability measure $\mu$ on $\R^n$ is called isotropic if for
   every $y \in \R^n$, $\E |\inr{X,y}|^2 = |y|^2$, where $X$ is
   distributed according to $\mu$.

   A measure $\mu$ satisfies a
   $\psi_2$ condition with a constant $\alpha$ if for every $y \in
   \R^n$, $\|\inr{X,y}\|_{\psi_2} \leq \alpha |y|$.

   A subgaussian
   or $\psi_2$ operator is a random operator of the form
   \begin{equation}\label{Gamma}
   \Gamma=\sum_{i=1}^k \inr{X_i,\,.\,}e_i
   \end{equation}
    where the $X_i$
   are distributed according to an isotropic $\psi_2$ measure.
\end{Definition}

Perhaps the most important example of an isotropic $\psi_2$
probability measure on $\R^n$ with a bounded constant other than the
Gaussian measure is the uniform measure on $\{-1,1\}^n$.  Naturally,
if $X$ is distributed according to a general isotropic $\psi_2$
measure then the coordinates of $X$ need no longer be independent.
For example, the normalized Lebesgue measure on an appropriate
multiple of the unit ball in $\ell_p^n$ for $2 \leq p \leq \infty$ is
an isotropic $\psi_2$ measure with a constant independent of $n$ and
$p$.  For more details on such measures see \cite{MP}.

For a set $T \subset \R^n$ and $\rho >0$ let
\begin{equation}
  \label{Trho}
      T_\rho = T \cap \rho S^{n-1}.
\end{equation}

The next result shows that given $T\subset \R^n$, subgaussian operators
are very close to being an isometry on the subset of elements of $T$ which
have a ``large enough" norm.

\begin{Theorem}
   \label{thm:psi2}
   There exist absolute constants  $c, \bar{c} >0$
   for which the following holds.  Let $T \subset \R^n$ be a
   star-shaped set.  Let $\mu$ be an isotropic $\psi_2$ probability
   measure with constant $\alpha \ge 1$.  Let $ k \ge 1$, and
   $X_1, \ldots, X_k$ be independent, distributed according to $\mu$
   and define $\Gamma$ by (\ref{Gamma}). For   $0 < \theta <1$,
   with probability at least $1- \exp(- \bar{c}\, \theta^2 k/\alpha ^4)$,
   then  for all $x \in T$ such that  $|x| \ge r_k^*(\theta)$, we have
\begin{equation}
  \label{two_sided}
   (1 - \theta) |x|^2 \le \frac{|\Gamma x|^2}{ k} \le (1+\theta) |x|^2,
\end{equation}
where
\begin{equation}
  \label{rho_star}
r^*_k(\theta) = r^*_k (\theta, T) := \inf \Bigl\{ \rho>0 : \rho \ge
c \, \alpha^2  \lstar (T_\rho)/ \bigl(\theta \sqrt k \bigr) \Bigr\}.
\end{equation}
In particular, with the same probability, every
$x \in T$ satisfies
$$
|x|^2 \leq \max \Bigl\{(1-\theta)^{-1} |\Gamma x|^2/k, \
         r^*_k(\theta)^2 \Bigr\}.
$$
\end{Theorem}

\proof We use Theorem~\ref{thm:general} for the set of functions $F $
consisting of linear functionals of the form $f = f_x = \inr{\cdot,
  x}$, for $x \in T$.  By the isotropicity of $\mu$, $\|f\|_{L_2} =
|x| $ for $f = f_x \in F $.  Also, since $\mu$ is
$\psi_2$ with constant $\alpha$ then it follows by (\ref{maj_meas})
that for all $\rho >0$,
$$
\gamma_2(F \cap \rho S_{L_2},
\|\ \|_{\psi_2}) \le \alpha\, \gamma_2(F \cap \rho S_{L_2}, \|\ \|_{L_2})
\le ( \alpha /c_2) \, \lstar (T_\rho),
$$
as promised.
\endproof

\begin{Remark}
  \label{theta-large}
  {\rm It is clear from the proof of Theorem~\ref{thm:general} that
    the upper estimates in (\ref{two_sided}) hold for $\theta \ge 1$
    as well, with appropriate probability estimates and a modified
    expression for $r_k^*$ in (\ref{rho_star}). Note that of course in
    this case the lower estimate in (\ref{two_sided}) became vacuous.
    The same remark is valid for the estimate in
    (\ref{two_sided_gen}) as well. }
\end{Remark}

The last result immediately leads to an  estimate for the diameter of
random sections of a set  $T$ in $\R^n$, given by kernels of random
operators $\Gamma$, and which is a $\psi_2$-counterpart of the main
result from \cite{PT} (see also \cite{PTsem}).

\begin{Corollary}
  \label{sections}
   There exist absolute constants  $\tilde{c}, \tilde{c}' >0$
   for which the following holds.  Let $T \subset \R^n$ be a
   star-shaped set and  let $\mu$ be an isotropic $\psi_2$ probability
   measure with constant $\alpha \ge 1$.  Set $ k \ge 1$, put
   $X_1, \ldots, X_k$ to be independent, distributed according to $\mu$
   and define $\Gamma$ by (\ref{Gamma}).  Then, with probability at least $1
   - \exp(- \tilde{c} k /\alpha^4)$,
$$
\diam(\ker \Gamma \cap T) \le r_k^*(T),
$$
where
\begin{equation}
  \label{rho_star2}
r^*_k = r^*_k ( T) := \inf \Bigl\{ \rho>0 : \rho \ge
\tilde{c}' \, \alpha^2  \lstar (T_\rho)/ \sqrt k  \Bigr\}.
\end{equation}
Moreover, with the same probability,
$\diam(\ker \Gamma \cap T) \le
\tilde{c}' \, \alpha^2  \lstar (T)/ \sqrt k $.
\end{Corollary}

The Gaussian case (that is, when $\mu$ is the standard Gaussian
measure on $\R^n$), although not explicitly stated in \cite{PT},
follows immediately from the proof in that paper.  The parameter
$\inf\{\rho>0 : \lstar(T\cap \rho B_2^n)\le C\rho\sqrt k \}$ was
introduced in \cite{PTsem}.

A version of Corollary~\ref{sections} for random $\pm1$-vectors
follows from the result in \cite{Art}, as observed in \cite{MP}.

\medskip

\noindent{\bf Proof of Corollary~\ref{sections}.}
Applying  Theorem~\ref{thm:psi2} with $\theta = 1/2$, say, we get that
if $x \in T$ and $|x| \ge r_k^*(1/2)$ then $\Gamma x \ne 0$.  Thus for
$ x \in \ker \Gamma \cap T$ we have $|x| \le  r_k^*(1/2)$ and the first
conclusion follows by adjusting the constants.

Observe that since the function $\lstar(T_\rho)/\rho =
\lstar\bigl((1/\rho)T \cap S^{n-1}\bigr)$ is decreasing in $\rho$ then
$r_k^*$ actually satisfies the equality in the defining formula
(\ref{rho_star2}). Combining this and the obvious estimate $\lstar
(T_\rho) \le \lstar (T)$,  concludes  the ``moreover'' part.
\endproof

Finally, let us note a special case of Theorem~\ref{thm:psi2} for
subsets of the sphere.

\begin{Corollary}
  \label{expect}
  Let $T \subset S^{n-1}$ and let $\mu$, $\alpha$, $k$, $X_i$,
  $\Gamma$ and $\theta$ be the same as in Theorem~\ref{thm:psi2}.  As
  long as $k$ satisfies $k \ge (c' \, \alpha^4/\theta^2) \lstar (T
  )^2$, then with probability at least $1- \exp(- \bar{c}\, \theta^2
  k/\alpha ^4)$, for all $x \in T$,
\begin{equation}
  \label{two_sided-2}
   1 - \theta \le \frac{|\Gamma x|}{ \sqrt k} \le 1+\theta,
\end{equation}
where $c, \bar{c} >0$ are absolute constants.
\end{Corollary}

\proof Let $c, \bar{c} >0$ be the constants from
Theorem~\ref{thm:psi2}. Observe that the condition on $k$, with $c' =
c^2$, is equivalent to $r_k^*(\theta, \tilde{T}) \le 1$, where
$\tilde{T} = \{ \lambda x : x \in T, 0 \le \lambda \le 1 \}$.
Then (\ref{two_sided-2}) immediately follows from (\ref{two_sided}).
\endproof

\section{Approximate reconstruction}
\label{sec:app-rec}

Next, we show how one can apply Theorem \ref{thm:psi2} to
reconstruct any fixed $v \in T$ for any set $T \subset \R^n$, where
the data at hand are linear subgaussian measurements of the form
$\inr{X_i,v}$.

The reconstruction algorithm we choose is as follows: for a fixed
$\eps>0$, find some $t \in T$ such that
$$
\left(\frac{1}{k}\sum_{i=1}^k \left(\inr{X_i,v} -
\inr{X_i,t}\right)^2\right)^{1/2} \leq \eps.
$$
The fact that we only need to find $t$ for which
$(\inr{X_i,t})_{i=1}^k$ is close to
$(\inr{X_i,v})_{i=1}^k$ rather than equal to it, is very important
algorithmically  because it is a far simpler problem.

Let us show why such an algorithm can be used to solve the
approximate reconstruction problem.

Consider ${\bar T}=\left\{\lambda(t-s):t,s \in T, \ 0 \leq \lambda
  \leq 1\right\}$ and observe that by Theorem \ref{thm:psi2}, for every
$0 < \theta <1$, with high probability, every such $t \in T$ satisfies
that
$$
|t-v| \leq \frac{\eps}{1-\theta} + r_k^*(\theta,\bar{T}).
$$
Hence, to bound the reconstruction error, one needs to estimate
$r_k^*(\theta,\bar{T})$. Of course, if $T$ happens to be convex and
symmetric then ${\bar T} \subset 2T$ which is star-shaped and thus
$$
|t-v| \leq \frac{\eps}{1-\theta}+r^*_k(\theta,2T).
$$
In a more general case, when $T$ is symmetric and quasi-convex with
constant $a \ge 1$, (i.e., $T + T \subset 2 a T$  and $T$ is
star-shaped), then
$$
|t-v| \leq \frac{\eps}{1-\theta} + r_k^*(\theta, a T).
$$

Therefore, in the quasi-convex case, the ability to approximate any
point in $T$ using this kind of random sampling depends on the
expectation of the supremum of a Gaussian process indexed by the
intersection of $T$ and a sphere of a radius $\rho$ as a function of
the radius. For a general set $T$, the reconstruction error is
controlled by the behavior of the expectation of the supremum of the
Gaussian process indexed by the intersection of $\bar{T}$ with
spheres of radius $\rho$, and this function of $\rho$ is just the
modulus of continuity of the Gaussian process indexed by the set
$\{\lambda t : 0 \leq \lambda \leq 1, \ t \in T \}$ (i.e, the
expectation of the supremum of the Gaussian process indexed by the
set $\{\lambda(t-s) : 0 \leq \lambda \leq 1, \ t,s \in T, \ |t-s| =
\rho\}$).

The parameters $r_k^*(\theta, T)$ can be estimated for the unit ball
of classical normed or quasi-normed spaces.  The
two example we consider here are the unit ball in $\ell_1^n$, denoted
by $B_1^n$, and the unit balls in the weak-$\ell_p^n$ spaces
$\ell_{p,\infty}^n$ for $0<p<1$, denoted by $B_{p,\infty}^n$.  Recall
that $B_{p,\infty}^n$ is the set of all $x = (x_i)_{i=1}^n \in \R^n$
such that the cardinality $|\{i: |x_i|\ge s\}|\le s^{-p}$ for all $
s>0$, and observe that $B_{p,\infty}^n$ is a quasi convex body with
constant  $a = 2^{1/p}$. Let us mention that there is nothing
``magical" about the examples we consider here. Those are simply the
cases considered in \cite{CT1,RV}.

In order to bound $r_k^*$ for these sets we shall use the approach
from \cite{GLMP}, and combine it  with Theorem \ref{thm:psi2} to
recover and extend the results from \cite{CT1,RV}.

\begin{Theorem}
  \label{cor:B_1^n}
  There is an absolute constant ${\bar c}$ for which the following
  holds.  Let $1 \le k \le n$ and $0<\theta <1$, and set $\eps>0$.
  Let $\mu$ be an isotropic $\psi_2$ probability measure on $\R^n$
  with constant $\alpha$, and let $X_1, \ldots, X_k$ be independent,
  distributed according to $\mu$. For any $0 < p <1 $, with
  probability at least $1-\exp(-{\bar{c}}\theta^2k/\alpha^4)$, if $v,y
  \in B_{p, \infty}^n$ satisfy
  that $(\sum_{i=1}^k \inr{X_i,v-y}^2 /k)^{1/2} \leq
  \eps$, then
$$
|y-v| \leq   \frac{\eps}{1-\theta} +
2^{1/p +1} \Bigl(\frac{1}{p} -1\Bigr)^{-1}
\left(C_{\alpha, \theta}\, \frac{\log(C_{\alpha, \theta}n/k)}{k}\right)
^{1/p - 1/2},
$$
where $C_{\alpha, \theta} = c \alpha^4/\theta^2$ and $c >0$ is an
absolute constant.

If $v,y \in B_1^n$ satisfy the same assumption then with the same
probability estimate,
$$
|y-v| \leq \frac{\eps}{1-\theta} +
 \left( C_{\alpha, \theta}\,\frac{\log(C_{\alpha, \theta}n/k)}{k}
\right)^{1/2}.
$$
\end{Theorem}

To prove Theorem \ref{cor:B_1^n}, we require the following
elementary fact.
\begin{Lemma}
   \label{lemma:inter}
   Let $0 < p <1$ and $1 \leq m \leq n$.  Then, for every $x \in
   \R^n$,
$$
\sup_{z  \in B_{p,\infty}^n\cap \rho B_2^n} \inr{x,z } \leq 2\rho
\left(\sum_{i=1}^m {x_i^*}^2 \right)^{1/2},
$$
where
   $\rho=\left(1/p-1\right)^{-1}m^{1/2-1/p}$ and
$(x_i^*)_{i=1}^n$ is a non-increasing rearrangement of
$(|x_i|)_{i=1}^n$.

Moreover,
$$
\sup_{z  \in B_1^n \cap \rho B_2^n } \inr{x,z } \leq 2\rho
\left(\sum_{i=1}^m {x_i^*}^2 \right)^{1/2},
$$
with $\rho=1/\sqrt{m}$.
\end{Lemma}

\proof We will present a proof for the case $0<p<1$. The case of
$B_1^n$ is similar.

Recall a well known fact that for two sequences of positive numbers
$a_i, b_i$ such that $a_1 \ge a_2 \ge \ldots$, the sum
$\sum a_i b_{\pi(i)}$ is maximal over all permutations $\pi$ of the
index set, if
$b_{\pi(1)} \ge b_{\pi(2)} \ge \ldots$. It follows that,
 for any $\rho>0, m\ge 1$ and $z  \in
B_{p,\infty}^n \cap \rho B_2^n$,
\begin{align*}
\inr{x,z }  \leq & \rho \left(\sum_{i=1}^m {x_i^*}^2
\right)^{1/2} +\sum_{i>m} \frac{x_i^*}{i^{1/p}} \\
\leq & \left(\sum_{i=1}^m {x_i^*}^2 \right)^{1/2}
\left(\rho+\frac{1}{\sqrt{m}}\sum_{i >m} \frac{1}{i^{1/p}} \right)
\\
\leq & \left(\sum_{i=1}^m {x_i^*}^2 \right)^{1/2}
\left(\rho+\left(\frac{1}{p}-1\right)^{-1}
\frac{1}{m^{1/p-1/2}}\right).
\end{align*}
By the definition of $\rho$, this completes the proof.
\endproof

\medskip

Consider the set of elements in the unit ball with ``short support",
defined by
$$
U_m=\left\{x \in S^{n-1} : \left| \left\{i: x_i \not = 0 \right\}
\right| \leq m \right\}.
$$

\smallskip

Note that Lemma \ref{lemma:inter} combined with a duality argument
implies that for every $1 \leq m \leq n$ and every $I \subset
\{1,...,n\}$ with $|I| \leq m$,
\begin{equation} \label{eq:Um}
\sqrt{|I|}\, B^n_1\cap S^{n-1}\subset 2 \conv  U_{m} \cap
S^{n-1}.
\end{equation}

\smallskip

The next step is to bound the expectation of the supremum of the
Gaussian process indexed by $U_m$.

\begin{Lemma}
   \label{lemma:gauss-Um}
   There exist an absolute constant $c$ such that for every $1 \leq m
   \leq n$,
$$
\lstar(\conv U_m)  \leq c \sqrt{m\log(cn/m)}.
$$
\end{Lemma}

\proof Recall that for every $1 \leq m \leq n$, there is a set
$\Lambda_m$ of cardinality at most $5^m$ such that $B_2^m \subset 2
\conv \Lambda_m$ (for example, a successive approximation shows that
we may take as $\Lambda_m$ an $1/2$-net in $B_2^m $).
Hence there is a subset of $B_2^n$ of cardinality at most $5^m
\binom{n}{m}$ such that $U_m \subset 2\conv \Lambda $. It is well
known (see for example\cite{LT}) that for every $T \subset B_2^n$,
$$
\lstar(\conv T)=\lstar (T) \leq c\sqrt{\log(|T|)},
$$
and thus,
$$
\lstar (\conv U_m) \leq c\sqrt{\log\Bigl(5^m \binom{n}{m}\Bigr)},
$$
from which the claim follows.
\endproof

Finally, we are ready to estimate $r_k^*(\theta,B_{p,\infty}^n)$ and
$r_k^*(\theta,B_1^n)$.

\begin{Lemma}
   \label{lemma:r-est}
There exists an absolute constant $c$ such that for any $0<p<1$ and
$1 \leq k \leq n$,
\begin{equation*}
r_k^*(\theta,B_{p,\infty}^n) \leq c\left(\frac{1}{p}-1\right)^{-1}
\left(\frac{\log(cn\alpha^4/\theta^2k)}{\theta^2k/\alpha^4}\right)^{1/p-1/2}
\end{equation*}
and
\begin{equation*}
r_k^*(\theta,B_1^n) \leq c
\left(\frac{\log(cn\alpha^4/\theta^2k)}{\theta^2k/\alpha^4}\right)^{1/2}.
\end{equation*}
\end{Lemma}

\proof Again, we present a proof for $B_{p,\infty}^n$, while the
treatment of $B_1^n$ is similar and thus omitted.

Let $0<p<1$ and $1 \leq k \leq n$, and set $1 \leq m \leq n$ to be
determined later. Clearly,
$$
\left(\sum_{i=1}^m {x_i^*}^2 \right)^{1/2} = \sup_{y \in U_m}
\inr{x,y},
$$
and
thus, by Lemma \ref{lemma:inter}, $\lstar(B_{p,\infty}^n \cap \rho
B_2^n) \leq 2\rho\lstar(U_m)$, where
$\rho=(1/p-1)^{-1}m^{1/2-1/p}$. From the definition of $r_k^*(\theta)$ in
Theorem \ref{thm:psi2}, it suffices to determine $m$ (and thus $\rho$)
such that
$$
c\alpha^2\lstar(U_m) \leq \theta \sqrt{k},
$$
which by Lemma \ref{lemma:gauss-Um}, comes to $c\alpha^2
\sqrt{m\log(cn/m)} \leq \theta \sqrt{k}$ for some other numerical
constant $c$. It is standard to verify that the last inequality
holds true provided that
$$
m \leq
c\frac{\theta^2k/\alpha^4}{\log\left(cn\alpha^4/\theta^2k\right)},
$$
and thus
$$
r_k^*(\theta,B_{p,\infty}^n) \leq c \left(\frac{1}{p}-1\right)^{-1}
\left(\frac{\log\left(cn\alpha^4/\theta^2k\right)}
   {\theta^2k/\alpha^4}\right)^{1/p-1/2}.
$$
\endproof

\medskip
\noindent {\bf Proof of Corollary \ref{cor:B_1^n}.} The proof
follows immediately from Theorem \ref{thm:psi2} and Lemma
\ref{lemma:r-est}.
\endproof

\section{Exact reconstruction}
\label{sec:exact}

Let us consider the following problem from the error correcting code
theory. A linear code is given by an $n\times (n-k)$ real matrix
$A$. Thus, a vector $x\in \R^{n-k}$ generates the vector
$Ax\in\R^n$. Suppose that $Ax$ is corrupted by a noise vector
$z\in\R^n$ and the assumption we make is that $z$ is {\em sparse},
that is, has a short support, which we denote by supp$(z)$=$\{i\,:\,
z_i\not=0\}$. The problem is to {\em reconstruct} $x$ from the data,
which is the noisy output $y=Ax+z$.

For this purpose, consider a $k\times n$ matrix $\Gamma$ such that
$\Gamma A=0$. Thus $\Gamma z=\Gamma y$ and correcting the noise is
reduced to identifying the sparse vector $z$ (rather than
approximating it) from the data $\Gamma z$ - which is the problem we
focus on here.

In this context, a linear programming approach called the {\em basis
  pursuit} algorithm, was recently shown to be relevant
for this goal \cite{CDS}. This method is the following minimization
problem
$$
(P)\qquad  \min \|t\|_{\ell_1},\quad \Gamma t=\Gamma z
$$
(and recall that the ${\ell_1}$-norm is defined by
$\|t\|_{\ell_1}=\sum_{i=1}^n |t_i|$ for any $t=(t_i)_{i=1}^n
\in\R^n$).

For the analysis of the reconstruction of sparse vectors by this basis
pursuit algorithm, we refer to \cite{CDS} and the recent papers
\cite{CT2, CT3}.

In this section, we show that if $\Gamma$ is an isotropic $\psi_2$
matrix then with high probability, for any vector $z$ whose support
has size less than ${Ck\over \log (Cn/k)}$ (for some absolute
constant $C$), the problem $(P)$ above has a unique solution that is
equal to $z$. It means that such random matrices can be used to
reconstruct any sparse vector, as long as the size of the support is
not too large. This extends the recent result proved in \cite{CT2}
and \cite{RV} for Gaussian matrices.

\begin{Theorem} \label{thm:reconstruction} There exist absolute
  constants $c,C$ and $\bar{c}$ for which the following holds. Let
  $\mu$ be an isotropic $\psi_2$ probability measure with constant
  $\alpha \ge 1$.  For  $ 1\le k \le n$, set $X_1, \ldots, X_k$ to be
  independent, distributed according to $\mu$ and let
  $\Gamma=\sum_{i=1}^k \inr{X_i,\cdot}e_i$. Then with probability at
  least $1- \exp(- \bar{c}\, k/\alpha ^4)$, any vector $z$ satisfying
$$
|\supp(z)|\le {C k\over \alpha^4 \log (cn \alpha^4/k)}
$$
is the unique minimizer of the problem
$$
(P)\qquad  \min \|t\|_{\ell_1},\quad \Gamma t=\Gamma z.
$$
\end{Theorem}

The proof of Theorem \ref{thm:reconstruction} is based on
a scheme of
the proof of \cite{CT2}, but is simpler and more general, as it
holds for an arbitrary isotropic $\psi_2$ random matrix.

Let us remark that problem $(P)$ is equivalent to the following one
$$(P')\qquad
\min_{t\in \R^n} \|y-At\|_{\ell_1}
$$
where $\Gamma A=0$. Thus we obtain the reconstruction result:

\begin{Corollary}\label{thm:reconstruction2}
Let $A$ be a $n\times (n-k)$ matrix. Set $\Gamma$ to be a $k\times
n$ matrix that satisfies the conclusion of the previous Theorem
with the constants $c$ and $C$, and for which $\Gamma A=0$. For
any $x\in \R^{n-k}$ and any $y=Ax+z$, if $|{\rm supp}(z)|\le
{Ck\over \log (cn/k)}$, then $x$ is the unique minimizer of the
problem
$$
\min_{t\in \R^n} \|y-At\|_{\ell_1}.
$$
\end{Corollary}

The condition $\Gamma A=0$ means that the range of $A$ is a subspace
of the kernel of $\Gamma$. Due to the rotation invariance of a
Gaussian matrix (from both sides), the range and the kernel are
random elements of the Grassmann manifold of the corresponding
dimensions. Therefore, random Gaussian matrices $A$ and $\Gamma$
satisfy the conclusion of Corollary~\ref{thm:reconstruction2}.

\subsection{Proof of Theorem \ref{thm:reconstruction}}

As in \cite{CT2}, the proof consists of finding a simple condition
for a fixed matrix $\Gamma$ to satisfy the conclusion of our
Theorem. We then apply a result from the previous section to show
that random matrices satisfy this condition.

The first step is to provide some criteria which ensure that the
problem $(P)$ has a unique solution as specified in Theorem
\ref{thm:reconstruction}. This convex optimization problem can be
represented as a linear programming problem. Indeed, let $z\in\R^n$
and set
\begin{equation}
  \label{supports}
I^+=\{i\,:\ z_i>0\},\quad I^-=\{i\,:\ z_i<0\}, \quad I=I^+\cup I^- .
\end{equation}
Denote by  $ {\cal C}$ the cone of constraint
$$ {\cal C}=\{t\in\R^n\,:\, \sum_{i\in I^+} t_i-\sum_{i\in I^-} t_i+
\sum_{i\in I^c} |t_i|\le 0\}$$ corresponding to the $\ell_1$-norm.

Note that if $|t|$ is small enough
then
$\|z+t\|_{\ell_1}=
\sum_{i\in I^+}(z_i + t_i)-\sum_{i\in I^-} (z_i+ t_i)+
\sum_{i\in I^c} |t_i|$.
Thus, the solution of $(P)$ is unique and
equals to $z$ if and only if
\begin{equation}
\label{constraint condition}
\ker \Gamma\cap {\cal C}=\{0\}
\end{equation}

By the Hahn-Banach separation Theorem, the latter is equivalent to
the existence of a linear form $\tilde w\in \R^n$ vanishing on
$\ker \Gamma$ and positive on ${\cal C}\setminus\{0\}$.

After appropriate normalization, it is easy to check that such an
$\tilde w$ satisfies that $\tilde w=\sum_{i=1}^k\alpha_i X_i$,
$\inr{\tilde w,e_i}=1$ for all $i\in I^+$, $\inr{\tilde w,e_i}=-1$ for
all $i\in I^-$, and $|\inr{\tilde w,e_i}|<1$ for all $i\in I^c$.
Setting $w=\sum_{i=1}^k\alpha_i e_i$ and noticing that $\inr{\tilde
  w,e_i}=\inr{w,\Gamma e_i}$ we arrive at the following criterion.

\begin{Lemma}
\label{face condition}
  Let $\Gamma$ be a $k\times n$ matrix and $z\in\R^n$. With the
  notation (\ref{supports}), the problem
 $$
(P)\qquad  \min \|t\|_{\ell_1},\quad \Gamma t=\Gamma z
$$
has a unique solution which equals to $z$, if and only if there
exists $w\in\R^k$ such that
$$
\forall i\in I^+ \ \inr{w,\Gamma e_i}=1,\quad
\forall i\in I^- \ \inr{w,\Gamma e_i}=-1,\quad \forall i\in I^c \
|\inr{ w,\Gamma e_i}|<1.
$$
\end{Lemma}

The second preliminary result we require follows from Corollary
\ref{expect} and the estimates of the previous section.

\begin{Theorem}
  \label{thm:Um}
  There exist absolute constants $c, C$ and $\bar{c}$ for which the
  following holds. Let $\mu$, $\alpha$, $k$ and $\Gamma$ be as in
  Theorem~\ref{thm:reconstruction}.  Then, for every $0 < \theta <1$,
  with probability at least $1- \exp(- \bar{c}\, \theta^2 k/\alpha
  ^4)$, every $x \in 2 \conv U_{4m} \cap S^{n-1}$ satisfies that
\begin{equation}
  \label{lower}
   (1 - \theta) |x|^2 \le \frac{|\Gamma x|^2}{ k} \le (1+\theta) |x|^2,
\end{equation}
provided that
$$
m \leq
C\frac{\theta^2k/\alpha^4}{\log\left(cn\alpha^4/\theta^2k\right)}.
$$
\end{Theorem}

\noindent {\bf Proof.} Applying Corollary \ref{expect} to $T=2 \conv
U_{4m} \cap S^{n-1}$, we only have to check that $k \ge (c' \,
\alpha^4/\theta^2) \lstar (T )^2$, which from Lemma
\ref{lemma:gauss-Um} reduces to verifying that $ k\ge (c'\alpha^4
/\theta^2)cm\log(cn/m) $. The conclusion now follows from the same
computation as in the proof of Lemma \ref{lemma:r-est}.

\bigskip \noindent {\bf Proof of Theorem \ref{thm:reconstruction}.}
Observe that if $t\in {\cal C}\cap S^{n-1}$ then
$\|t\|_{\ell_1}\le 2\sum_{i\in I} |t_i|\le 2\sqrt{|I|}$, where $I$
is the support of $z$. Hence,
$$
{\cal C}\cap S^{n-1}\subset \sqrt{4|I|}\, B^n_1\cap S^{n-1}.
$$
This inclusion and condition (\ref{constraint condition}) clearly
imply that if $\Gamma$ does not vanish on any point of $\sqrt{4|I|}\,
B^n_1\cap S^{n-1}$, then the solution of $(P)$ is unique and equals to
$z$. By \eqref{eq:Um} we have
$$
\sqrt{4|I|}\, B^n_1\cap S^{n-1}\subset 2\conv U_{4m} \cap
S^{n-1}.
$$
Therefore, if $\Gamma$ does not vanish on any point of $2\conv
U_{4m}\cap S^{n-1}$ then $z$ is the unique solution of $(P)$.
Applying Theorem \ref{thm:Um}, the lower bound in (\ref{lower})
shows that indeed, $\Gamma$ does not vanish on any point of the
required set, provided that
$$
m\le  {Ck\over \alpha^4\log (cn \alpha^4/k)}
$$
for some suitable constants $c$ and $C$.
\endproof

\subsection{The geometry of faces of random polytopes}
\label{faces}

Next, we investigate the geometry of random polytopes. Let $\Gamma$ be
a $k\times n$ isotropic $\psi_2$ matrix. For $1 \leq i \leq n$ let
$v_i=\Gamma(e_i)$ be the vector columns of the matrix $\Gamma$ and set
$K^+(\Gamma)$ (resp. $K(\Gamma)$) to be the convex hull (resp., the
symmetric convex hull) of these vectors.

In this situation, the random model that makes sense is when
$X=(x_i)_{i=1}^n$, where $(x_i)_{i=1}^n$ are independent,
identically distributed random variables for which $\E|x_i|^2 =1$
and $\|x_i\|_{\psi_2} \leq \alpha$. It is standard to verify that in
this case $X =(x_i)_{i=1}^n $ is an isotropic $\psi_2$ vector with constant
$\alpha$, and moreover, each vertex of the polytope is given by
$v_i=(x_{i,j})_{j=1}^k$.

A polytope is called {\em $m$-neighborly} if any set of less than $m$
vertices is  the vertex set of a face.
In the symmetric setting, we will say that 
a symmetric polytope is
$m$-symmetric-neighborly if any set of less than $m$
vertices containing no-opposite pairs,
is  the vertex set of a face.

The condition of Lemma \ref{face condition} may be reformulated by
saying that the set $\{v_i\,:\, i\in I^+\}\cup \{-v_i\,:\, i\in
I^-\}$ is the vertex set of a face of the polytope $K(\Gamma)$.
Thus, the condition for the exact reconstruction using the basis
pursuit method for any vector $z$ with  $|\supp (z)| \leq m$ may
be reformulated as a geometric property of the polytope
$K(\Gamma)$ (see \cite{CT2,RV}); namely, that {\em for all
disjoint subsets $I^+$ and $I^-$ of $\{1,\dots, n\}$ such that
$|I^+|+|I^-|\le m$, the set $\{v_i\,:\, i\in I^+\}\cup \{-v_i\,:\,
i\in I^-\}$ is the vertex set of a face of the polytope
$K(\Gamma)$.} That is, $K(\Gamma)$ is $m$-symmetric-neighborly. A similar
analysis may be done in the non-symmetric case, for $K^+(\Gamma)$,
where now $I^-$ is empty.

\begin{Lemma}
Let $\Gamma$, $K(\Gamma)$ and $K^+(\Gamma)$ be as above. Then the
problem
$$
(P)\qquad  \min \|t\|_{\ell_1},\quad \Gamma t=\Gamma z
$$
has a unique solution which equals to $z$ for any vector $z$ (resp.,
$z\ge 0$) such that $|{\rm\, supp} (z)|\le m$, if and only if
$K(\Gamma)$ (resp., $K^+(\Gamma)$\/) is $m$-symmetric-neighborly
(resp., $m$-neighborly).
\end{Lemma}

Applying Theorem \ref{thm:reconstruction}, we obtain

\begin{Theorem}
  \label{thm:neighborly}
  There exist absolute constants $c,C$ and $\bar{c}$ for which the
  following holds. Let $\mu$ be an isotropic $\psi_2$ probability
  measure with constant $\alpha \ge 1$ and let $k$ and $\Gamma$ be as
  above. Then, with probability at least $1- \exp(- \bar{c}\, k/\alpha
  ^4)$, the polytopes $K^+(\Gamma)$ and $K(\Gamma)$ are $m$-neighborly
  and $m$-symmetric-neighborly, respectively, for every $m$ satisfying
$$
m\le {Ck\over \alpha^4\log (cn \alpha^4/k)}\cdot
$$
\end{Theorem}

The statement of Theorem \ref{thm:neighborly}
for $K(\Gamma)$ and for a Gaussian matrix $\Gamma$ is the main result
of \cite{RV}. However, a striking fact is that the same results holds for a
random $\{-1,1\}$-matrix. In such a case, $K^+(\Gamma)$ is the
convex hull of $n$ random vertices of the discrete cube $\{-1,+1\}^k$,
also known as a random $\{-1,1\}$-polytope.
With high probability, every $(m-1)$-dimensional face
of $K^+(\Gamma)$ is a simplex  and there are ${n\choose m}$ such
faces, for $m\le Ck/ \log (cn/k)$.

\medskip

\begin{Remark} {\rm Let us mention some related results about random
    $\{-1,1\}$-polytopes.  A result of \cite{BP} states that for such
    polytopes, the number of facets , which are the $k-1$-dimensional
    faces, may be super-exponential in the dimension $k$, for an
    appropriate choice of the number $n$ of vertices.  Denote by
    $f_{q}(K^+(\Gamma))$ the number of $q$-dimensional faces of the
    polytope $K^+(\Gamma)$. The quantitative estimate in \cite{BP} was
    recently improved in \cite{GGM} where it is shown that there are
    positive constants $a,b$ such that for $k^a\le n\le \exp(bk)$, one
    has $\E f_{k-1}(K^+(\Gamma))\ge (\ln n/\ln k^a)^{k/2}.$ For lower
    dimensional faces a threshold of $f_{q}(K^+(\Gamma))/ {n \choose
      q+1}$ was established in \cite{K}.}
\end{Remark}

\footnotesize
{
}

\noindent {\bf S. Mendelson} {\footnotesize Centre for Mathematics
and its Applications, The
Australian National University, Canberra, ACT 0200,
Australia \\} {\small\tt%
shahar.mendelson@anu.edu.au}  \\ [.05cm]

\noindent {\bf A. Pajor }{\footnotesize Laboratoire d'Analyse et
Math\'ematiques Appliqu\'ees, Universit\'e de Marne-la-Vall\'ee, 5
boulevard Descartes, Champs sur Marne, 77454 Marne-la-Vallee,
Cedex 2, France \\ }
{\small\tt%
    alain.pajor@univ-mlv.fr}\\ [.05cm]

\noindent {\bf N. Tomczak-Jaegermann} {\footnotesize
Department of Mathematical and Statistical Sciences,\\
University of Alberta,
Edmonton, Alberta, Canada T6G 2G1\\ }
{\small\tt%
  nicole@ellpspace.math.ualberta.ca}

\end{document}